\documentclass[12pt]{amsart}
\usepackage{amssymb, amsmath}

\newtheorem{lema}{Lemma}[section]
\newtheorem{teo}{Theorem}[section]

\newtheorem{coro}{Corollary}[section]
\newcommand{\mb}[1]{\mathbb{#1}}
\newcommand{\mc}[1]{\mathcal{#1}}
\DeclareMathOperator{\vol}{vol}

\begin{document}
\title{An estimate of free entropy and applications}
\author{Marius B. Stefan}
\address{UCLA Mathematics Department, Los Angeles, CA 90095-1555}
\email{stefan@math.ucla.edu}
 \abstract
We obtain an estimate of free entropy of generators in a type
$\mbox{\!I\!I}_1$-factor $\mc{M}$ which has a subfactor $\mc{N}$
of finite index with a subalgebra
$\mc{P}=\mc{P}_1\vee\mc{P}_2\subset\mc{N}$ where
$\mc{P}_1=\mc{R}_1'\cap\mc{P}$, $\mc{P}_2=\mc{R}_2'\cap\mc{P}$ are
diffuse, $\mc{R}_1,\mc{R}_2\subset\mc{P}$ are mutually commuting
hyperfinite subfactors, and an abelian subalgebra
$\mc{A}\subset\mc{N}$ such that the correspondence
$_\mc{P}L^2(\mc{N},\tau )_\mc{A}$ is $\mc{M}$-weakly contained in
a subcorrespondence $_\mc{P}H_\mc{A}$ of $_\mc{P}L^2(\mc{M},\tau
)_\mc{A}$, generated by $v$ vectors. The (modified) free entropy
dimension of any generating set of $\mc{M}$ is $\leq 2r+2v+4$,
where $r$ is the integer part of the index. As a consequence, the
interpolated free group subfactors of finite index do not have
regular non-prime subfactors or regular diffuse hyperfinite
subalgebras.
\endabstract
\maketitle

\section{Introduction}\label{s1}

D. Voiculescu defined (\cite{12}, \cite{13}) the original concepts
of free entropy and of (modified) free entropy dimension for
$m$-tuples of self-adjoint non-commutative random variables. Very
roughly, the free entropy $\chi ((x_i)_{1\leq i\leq m})$ is a
normalized limit of logarithms of volumes of sets of matricial
microstates (that is, $m$-tuples of matrices whose non-commutative
moments approximate those of $(x_i)_{1\leq i\leq m}$), while the
modified free entropy dimension $\delta_0 ((x_i)_{1\leq i\leq m})$
is in some sense an asymptotic Minkowski dimension of the sets of
matricial microstates. Then he proved (\cite{13}) that
$\delta_0((x_i)_{1\leq i\leq m})\leq 1$ if the von Neumann algebra
$\{(x_i)_{1\leq i\leq m}\}''$ has a regular diffuse hyperfinite
$*$-subalgebra (DHSA). Since the free group factors have
generators with $\delta_0>1$, this implied in particular the
absence of Cartan subalgebras in the free group factors, thus
answering in the negative the longstanding open question of
whether every separable $\mbox{\!I\!I}_1$-factor arises from a
measurable equivalence relation.

A. Connes introduced Kazhdan's property T from groups (\cite{23})
to the von Neumann algebras context in \cite{21} where he proved
that if $\Gamma$ is a countable discrete ICC group with property
T, then the fundamental group of the von Neumann algebra of
$\Gamma$ is countable. This remarkable rigidity result produced
the first examples (such as $\mc{L}(SL(3,\mb{Z}))$) of
$\mbox{\!I\!I}_1$-factors with fundamental group strictly smaller
than $\mb{R}^\times_+$. A. Connes and V. Jones (\cite{19}) defined
then property T for arbitrary von Neumann algebras in terms of
correspondences (\cite{20}, \cite{15}) and showed that the
countable discrete ICC group $\Gamma$ has property T if and only
if $\mc{L}(\Gamma)$ has property T. Correspondences play the role
of group representations and property T can be naturally defined
by following this analogy. Thus, the space of (equivalence classes
of) correspondences can be endowed with a topology through their
coefficients and property T simply means that the trivial
correspondence is isolated from the correspondences that do not
contain it.

The free group factors $\mc{L}(\mb{F}_n)$, $2\leq n<\infty$, were
the first examples of prime (i.e., non-isomorphic to tensor
products of) type $\mbox{\!I\!I}_1$-factors with separable
preduals. Their primeness (conjectured by S. Popa in \cite{24})
was proved by L. Ge (\cite{16}) via an estimate of free entropy.
Also with free entropy estimates and extending D. Voiculescu's
result about the absence of Cartan subalgebras, L. Ge (\cite{3})
and K. Dykema (\cite{2}) showed that the free group factors do not
have abelian subalgebras of multiplicity one and of finite
multiplicity, respectively.

We consider subcorrespondences $_\mc{P}H_\mc{A}$ of
$_\mc{P}L^2(\mc{M},\tau )_\mc{A}$, where $\mc{N}\subset\mc{M}$ is
an inclusion of type $\mbox{\!I\!I}_1$-factors with finite Jones
index (with integer part equal to $r$) and $\mc{P}$, $\mc{A}$ are
von Neumann subalgebras of $\mc{N}$. We assume moreover that
$\mc{A}$ is abelian and $\mc{P}=\mc{P}_1\vee\mc{P}_2$ where
$\mc{P}_1=\mc{R}_1'\cap\mc{P}$ and $\mc{P}_2=\mc{R}_2'\cap\mc{P}$
are both diffuse and $\mc{R}_1$, $\mc{R}_2$ are mutually commuting
hyperfinite subfactors. Then we prove (Theorem \ref{t1}, using
Lemma \ref{l1}) that if the correspondence $_\mc{P}L^2(\mc{N},\tau
)_\mc{A}$ is $\mc{M}$-weakly contained in $_\mc{P}H_\mc{A}$ and if
$_\mc{P}H_\mc{A}$ is spanned by $v$ vectors, then the (modified)
free entropy dimension of any generating set of $\mc{M}$ is $\leq
2r+2v+4$. For $\mc{M}=\mc{L}(\mb{F}_t)$ (an interpolated free
group factor from \cite{1}, \cite{7}), the free entropy dimension
estimate implies (Theorem \ref{t2}) the absence of regular
non-prime subfactors and of regular hyperfinite diffuse
subalgebras (DHSA) in the subfactors
$\mc{N}\subset\mc{L}(\mb{F}_t)$ ($1<t\leq\infty$) of finite Jones
index. In particular, the interpolated free group subfactors of
finite index are not crossed products of non-prime subfactors or
hyperfinite diffuse subalgebras by properly outer actions of
countable discrete groups (Corollary \ref{c1}). We mention that
the Haagerup approximation property (\cite{18}), primeness, and
absence of abelian subalgebras of finite multiplicity (and thus of
Cartan subalgebras) are known to be preserved (\cite{9},
\cite{10}) to the interpolated free group subfactors of finite
index.

\section{Notations}

We shall use $\mc{M}$, $\mc{N}$, $\mc{P}$, $\mc{A}$ etc. to denote
(finite) von Neumann algebras. In particular, we use $\mc{M}$ for
a type $\mbox{\!I\!I}_1$-factor and $\mc{N}$ for a subfactor of
$\mc{M}$. Let $\mc{P}$ be a finite von Neumann algebra, endowed
with a normal faithful tracial state
$\tau:\mc{P}\rightarrow\mb{C}$. For any projection $p\in\mc{P}$ we
denote by $\mc{P}_p=p\mc{P}p$ the corresponding reduced von
Neumann algebra of $\mc{P}$. Note that the functional
$\tau_p:\mc{P}_p\rightarrow\mb{C}$,
$\tau_p(y)=\frac{1}{\tau(p)}\tau(y)$ $\forall y\in\mc{P}_p$ is a
normal faithful tracial state on $\mc{P}_p$. The completion of
$\mc{P}$ with respect to the $2$-norm $||x||_2=
\tau(x^*x)^\frac{1}{2}$ $\forall x\in\mc{P}$ is a Hilbert space,
denoted $L^2(\mc{P},\tau)$. If $\mc{N}\subset\mc{M}$ is an
inclusion of type $\mbox{\!I\!I}_1$-factors, then the Jones index
$[\mc{M}:\mc{N}]$ is, by definition (\cite{5}), the dimension
$\mbox{dim}_\mc{N}L^2(\mc{M},\tau)$ of the left $\mc{N}$-module
$L^2(\mc{M},\tau)$. We mention that the dimension
$\mbox{dim}_\mc{N}H$ of an arbitrary left $\mc{N}$-module $H$ was
introduced by F. Murray and J. von Neumann (\cite{22}) as the {\it
coupling constant} of $H$. While $\mbox{dim}_\mc{N}H$ assumes all
possible values from $[0,\infty]$, the celebrated result of V.
Jones (\cite{5}) shows that necessarily $[\mc{M}:\mc{N}]\in\left
\{4\cos^2\frac{\pi}{n}:n\geq2\right\}\cup [4,\infty]$. Many
interesting von Neumann algebras arise from  representations of
discrete groups. For example, let $\Gamma$ be a discrete group and
denote by $(\delta_\gamma)_{\gamma\in\Gamma}$ the standard
orthonormal basis in $l^2( \Gamma)$. If
$\lambda:\Gamma\rightarrow\mc{B}(l^2(\Gamma))$, $\lambda
_{\gamma_1}\delta_{ \gamma_2}=\delta_{\gamma_1\gamma_2}$ $\forall
\gamma_1,\gamma_2\in\Gamma$ denotes the left regular
representation of $\Gamma$, then the (left) group von Neumann
algebra is $\mc{L}(\Gamma)=\lambda_\Gamma''$. Moreover, it is
easily seen that $\mc{L}(\Gamma)$ is a type
$\mbox{\!I\!I}_1$-factor if $\Gamma$ is an ICC group (that is, all
nontrivial conjugacy classes of $\Gamma$ are infinite). In
particular, the free group $\mb{F}_n$ on $n$ generators ($2\leq
n\leq\infty$) is an ICC group and thus one obtains the free group
factors $\mc{L}(\mb{F}_n)$ ($2\leq n\leq \infty$). The crossed
product construction is yet another way to obtain von Neumann
algebras from a von Neumann algebra $\mc{Q}$ and a discrete group
$\Gamma$ acting by $*$-automorphisms on $\mc{Q}$ (that is, there
exists a group homomorphism
$\alpha:\Gamma\rightarrow\mbox{Aut}(\mc{Q})$). Briefly, the
crossed product algebra $\mc{Q}\times_\alpha\Gamma$ is a kind of
maximal von Neumann algebra generated by (a copy of) $\mc{Q}$ and
(a copy of) $\Gamma$, subject to the commutation relations $\gamma
x\gamma^{-1}=\alpha_\gamma(x)$ $\forall \gamma\in\Gamma$ $\forall
x\in\mc{Q}$. The action
$\alpha:\Gamma\rightarrow\mbox{Aut}(\mc{Q})$ is called properly
outer if each automorphism $\alpha_\gamma$,
$\gamma\in\Gamma\setminus\{e\}$, is properly outer. An
automorphism $\beta\in\mbox{Aut}(\mc{Q})$ is properly outer if for
any $x\in\mc{Q}$, $xy=\beta(y)x$ $\forall y\in\mc{Q}$ implies
$x=0$. The action $\alpha$ is said to be ergodic if the
fixed-point subalgebra $\mc{Q}^\alpha
=\{x\in\mc{Q}:\alpha_\gamma(x)=x\,\forall\gamma\in\Gamma\}$ is
trivial. It is well-known that the crossed product
$\mc{Q}\times_\alpha\Gamma$ is a factor if the action $\alpha$ is
properly outer and if its restriction to the centre of $\mc{Q}$ is
ergodic. If $\mc{Q}\subset\mc{P}$ is an inclusion of von Neumann
algebras, then the normalizer $N_\mc{P}(\mc{Q})$ of $\mc{Q}$ in
$\mc{P}$ is defined by
$N_\mc{P}(\mc{Q})=\{u\in\mc{P}:uu^*=u^*u=1,u\mc{Q}u^*=\mc{Q}\}$.
The algebra $\mc{Q}$ is said to be regular in $\mc{P}$ if
$N_\mc{P}(\mc{Q})''=\mc{P}$. For example, $\mc{Q}$ is always
regular in $\mc{Q}\times_\alpha\Gamma$.

\subsection{Correspondences}

The notion of correspondence between two von Neumann algebras
(with separable preduals) $\mc{P}$ and $\mc{A}$ was introduced by
A. Connes (\cite{20}, \cite{15}, \cite{17}). Thus, a
correspondence between $\mc{P}$ and $\mc{A}$ is a pair of mutually
commuting normal unital $*$-representations of $\mc{P}$ and
$\mc{A}^o$ (the opposite algebra of $\mc{A}$) on the same
(separable) infinite dimensional Hilbert space $H$. Two
correspondences $_\mc{P}H_\mc{A}$ and $_\mc{P}H'_\mc{A}$ are
equivalent if there exists a $(\mc{P},\mc{A})$-bilinear isometry
from $H$ onto $H'$. We denote by $\widehat{_\mc{P}H_\mc{A}}$ the
class of $_\mc{P}H_\mc{A}$ under this equivalence relation and by
$\mbox{Corr}(\mc{P},\mc{A})$ the set of all classes of
correspondences between $\mc{P}$ and $\mc{A}$. Let
$\widehat{_\mc{P}H_\mc{A}} \in\mbox{Corr}(\mc{P},\mc{A})$,
$\epsilon >0$, and $F\subset\mc{P}$, $E\subset\mc{A}$,
$S=\{\xi_1,\ldots,\xi_v\}\subset H$ be finite subsets. Define
\begin{eqnarray}&&
U\left(\widehat{_\mc{P}H_\mc{A}},\epsilon,F,E,S\right)=\big\{\widehat{_\mc{P}K_\mc{A}}
\in\mbox{Corr}(\mc{P},\mc{A})\,:\,\exists\eta_1,\ldots,\eta_v\in
K\,\\\nonumber &&\hspace{.5 cm}|(b\xi_ka,\xi_l)
-(b\eta_ka,\eta_l)|<\epsilon\,\forall b\in F\,\forall a\in
E\,\forall 1 \leq k,l\leq v\big\}\,.
\end{eqnarray}
Define also
$V\left(\widehat{_\mc{P}H_\mc{A}},\epsilon,F,E,S\right)$ as the
set of all classes of correspondences
$\widehat{_\mc{P}K_\mc{A}}\in\mbox{Corr}(\mc{P}, \mc{A})$ with the
property that there exists a surjective isometry $u:H\rightarrow
K$ such that $||bu(\xi_l)a-u(b\xi_la)||<\epsilon$ for all $b\in
F$, $a\in E$, $1\leq l\leq v$. Then $\mbox{Corr}(\mc{P},\mc{A})$
becomes a topological space with the topology for which the sets
$U\left(\widehat{_\mc{P}H_\mc{A}},\epsilon,F,E,S\right)$ (or
equivalently, the sets
$V\left(\widehat{_\mc{P}H_\mc{A}},\epsilon,F,E,S\right)$) form a
basis of neighborhoods. One has also a notion of weak
subequivalence: if $_\mc{P}L_\mc{A}$ and $_\mc{P}L'_\mc{A}$ are
two correspondences between $\mc{P}$ and $\mc{A}$, we say that
$_\mc{P}L_\mc{A}$ is weakly contained in (or weakly subequivalent
to) $_\mc{P}L'_\mc{A}$ if
$\widehat{_\mc{P}L'_\mc{A}}\in\left\{\widehat{_\mc{P}L_\mc{A}}\right\}^-$
(we denoted by $V^-$ the closure of
$V\subset\mbox{Corr}(\mc{P},\mc{A})$).

We define next a refinement of the above equivalence relation,
restricted to subcorrespondences $_\mc{P}H_\mc{A}$ of
$_\mc{P}L^2(\mc{M},\tau )_\mc{A}$, where
$\mc{P}\vee\mc{A}\subset\mc{M}$. Thus, we say that two
subcorrespondences $_\mc{P}H_\mc{A}$, $_\mc{P}H'_\mc{A}$ of
$_\mc{P}L^2(\mc{M},\tau )_\mc{A}$ are $\mc{M}$-equivalent if there
exists a unitary $v\in(\mc{P}\vee\mc{A})'\cap\mc{M}$ such that
$H=vH'v^*$. We denote by $\widetilde{_\mc{P}H_\mc{A}}$ the class
of $_\mc{P}H_\mc{A}$ under this equivalence relation and by $
\mbox{Corr}_{\mc{M}}(\mc{P},\mc{A})$ the set of all classes
$\widetilde{_\mc{P}H_\mc{A}}$ of subcorrespondences
$_\mc{P}H_\mc{A}$ of $_\mc{P}L^2(\mc{M},\tau )_\mc{A}$. Denote
further by
$V_{\mc{M}}\left(\widetilde{_\mc{P}H_\mc{A}},\epsilon,F,E,S\right)$
the set of all classes $\widetilde{_\mc{P}K_\mc{A}}$ of
subcorrespondences $_\mc{P}K_\mc{A}$ of $_\mc{P}L^2(\mc{M},\tau
)_\mc{A}$ with the property that there exists a unitary
$w\in\mc{M}$ such that $K=wHw^*$ and
$||bw\xi_lw^*a-wb\xi_law^*||<\epsilon$ for all $b\in F$, $a\in E$,
$1\leq l\leq v$. Then $ \mbox{Corr}_{\mc{M}}(\mc{P},\mc{A})$
becomes a topological space, endowed with the topology for which
the sets
$V_{\mc{M}}\left(\widetilde{_\mc{P}H_\mc{A}},\epsilon,F,E,S\right)$
form a basis of neighborhoods. By analogy with the definition of
weak subequivalence, one defines the $\mc{M}$-weak subequivalence:
if $_\mc{P}L_\mc{A}$ and $_\mc{P}L'_\mc{A}$ are two
subcorrespondences of $_\mc{P}L^2(\mc{M},\tau )_\mc{A}$, we say
that $_\mc{P}L_\mc{A}$ is $\mc{M}$-weakly contained in (or
$\mc{M}$-weakly subequivalent to) $_\mc{P}L'_\mc{A}$ if
$\widetilde{_\mc{P}L'_\mc{A}}\in\left\{\widetilde{_\mc{P}L_\mc{A}}\right\}^-$
(where $V^-$ denotes here the closure of
$V\subset\mbox{Corr}_{\mc{M}}(\mc{P},\mc{A})$).

\subsection{Free entropy}

We recall but a few results from D. Voiculescu's free probability
theory (\cite{11}, \cite{12}, \cite{13}). Let $\mc{P}$ be a finite
von Neumann algebra endowed with a normal faithful tracial state
$\tau:\mc{P}\rightarrow\mb{C}$. An element $s\in\mc{P}$ is called
semicircular if $s=s^*$ and if it is distributed according to
Wigner's semicircle law:
\begin{equation}
\tau (s^k)=\frac{2}{\pi}\int_{-1}^{1}t^k\sqrt{1-t^2}dt\,\,\forall
k\in \mb{N}.
\end{equation}
A family of $(\mc{P}_i)_{i\in I}$ of unital $*$-subalgebras of
$\mc{P}$ is called free if $\tau(x_1\ldots$ $x_m)=0$ whenever
$x_k\in\mc{P}_{i_k}$, $\tau(x_k)=0$, $\forall 1\leq k\leq m$,
$i_1,\ldots,i_m\in I$, $i_1\not= i_2\not=\ldots\not= i_m$,
$m\in\mb{N}$. A family $(X_i)_{i\in I}$ of subsets
$X_i\subset\mc{P}$ is called free if the family
$(*$-$\mbox{alg}(\{1\}\cup X_i))_{i\in I}$ is free. The family
$(s_i)_{i\in I}$ of elements $s_i\in\mc{P}$ is a semicircular
system provided that $(\{s_i\})_{i\in I}$ is a free family and if
$s_i$ is a semicircular element $\forall i\in I$. If $c\geq 1$ is
an integer, we denote by $\mc{M}_c(\mb{C})$ and
$\mc{M}_c^{sa}(\mb{C})$ the set of all $c\times c$ complex
matrices and of all $c\times c$ complex self-adjoint matrices,
respectively. We further denote by $\mc{U}(c)$ the group of
unitaries from $\mc{M}_c(\mb{C})$, by $\tau_c$ the unique unital
trace on $\mc{M}_c(\mb{C})$, and by $||\cdot ||_e=\sqrt{c}||\cdot
||_2$ the euclidian norm on $\mc{M}_c(\mb{C})$. The free entropy
of $x_1,\ldots ,x_m\in\mc{P}^{sa}$ in the presence of
$x_{m+1},\ldots ,x_{m+n}\in\mc{P}^{sa}$ is defined in terms of
sets of matricial microstates $\Gamma_R((x_i)_{1\leq i\leq
m}:(x_{m+j})_{1\leq j\leq n};a,c,\epsilon )\subset
(\mc{M}_c^{sa}(\mb{C}))^m$. Thus, for $a,c\geq 1$ integers and
$R,\epsilon>0$, one has the following sequence of definitions:
\begin{eqnarray}&&
\Gamma_R\left((x_i)_{1\leq i\leq m}:(x_{m+j})_{1\leq j\leq
n};a,c,\epsilon\right)\\\nonumber &&\hspace{.5 cm}=
\left\{(A_i)_{1\leq i\leq m}\in (\mc{M}_c^{sa}
(\mb{C}))^m\,:\,\exists\,(A_{m+j})_{1\leq j\leq
n}\in(\mc{M}_c^{sa}(\mb{C}))^n\,\mbox{s.t.}\right.\\
\nonumber &&\hspace{.5 cm}|\tau \left.(x_{i_1}\ldots
x_{i_l})-\tau_c (A_{i_1}\ldots A_{i_l})|<\epsilon\,,\, ||A_k||\leq
R\,\right.\\\nonumber &&\hspace{.5 cm}\forall\left. 1\leq
i_1,\ldots ,i_l\leq m+n\,\forall 1\leq l\leq a\, \forall 1\leq
k\leq m+n\right\},
\end{eqnarray}
\begin{eqnarray}&&
\chi_R((x_i)_{1\leq i\leq m}:(x_{m+j})_{1\leq j\leq
n};a,c,\epsilon )\\\nonumber&&\hspace{.5 cm}
=\log\mbox{vol}_{mc^2}(\Gamma_R((x_i)_{1\leq i\leq
m}:(x_{m+j})_{1\leq j\leq n};a,c,\epsilon )),
\end{eqnarray}
\begin{eqnarray}&&\chi_R\left((x_i)_{1\leq i\leq m}:(x_{m+j})_{1\leq j\leq
n};a,\epsilon \right)\\\nonumber &&\hspace{.5 cm}
=\limsup\limits_{c\rightarrow\infty}\left(\frac{1}{c^2}\chi_R\left((x_i)_{1\leq
i\leq m}:(x_{m+j})_{1\leq j\leq n};a,c,\epsilon
\right)+\frac{m}{2}\log c\right),
\end{eqnarray}
\begin{eqnarray}&&
\chi_R\left((x_i)_{1\leq i\leq m}:(x_{m+j})_{1\leq j\leq
n}\right)\\\nonumber&&\hspace{.5 cm}
=\inf_{a,\epsilon}\chi_R\left((x_i)_{1\leq i\leq m}:
(x_{m+j})_{1\leq j\leq n};a,\epsilon\right),
\end{eqnarray}
\begin{eqnarray}&&\chi\left((x_i)_{1\leq i\leq m}:(x_{m+j})_{1\leq j\leq
n}\right)\\\nonumber&&\hspace{.5
cm}=\sup_R\chi_R\left((x_i)_{1\leq i\leq m}:(x_{m+j})_{1\leq j\leq
n} \right)
\end{eqnarray}
($\mbox{vol}_{mc^2}(\cdot )$ denotes the Lebesgue measure on
$(\mc{M}_c^{sa}(\mb{C}))^m\simeq\mb{R}^{mc^2}$). The last quantity
$\chi ((x_i)_{1\leq i\leq m}:(x_{m+j})_{1\leq j\leq n})$ is called
the free entropy of $(x_i)_{1\leq i\leq m}$ in the presence of
$(x_{m+j})_{1\leq j\leq n}$. If $n=0$, then it is simply called
the free entropy of $(x_i)_{1\leq i\leq m}$, denoted $\chi
(x_1,\ldots ,x_m)$. The free entropy of $(x_i)_{1\leq i\leq m}$ in
the presence of $(x_{m+j})_{1\leq j\leq n}$ is equal to the free
entropy of $(x_i)_{1\leq i\leq m}$ if $\{x_{m+1},\ldots
,x_{m+n}\}\subset\{x_1,\ldots ,x_m\}''$. For a single self-adjoint
element $x\in\mc{P}$ with distribution $\mu$ one has
\begin{equation}
\chi (x)=\frac{3}{4}+\frac{1}{2}\log 2\pi +\int\int\log |s-t|d\mu
(s)d\mu (t).
\end{equation}
Also, if $(x_i)_{1\leq i\leq m}$ is a free family, then $\chi
(x_1,\ldots ,x_m)=\chi (x_1)+\ldots +\chi (x_m)$. In particular, a
finite semicircular system has finite free entropy. The modified
free entropy dimension of $(x_i)_{1\leq i\leq m}$ is defined as
follows:
\begin{equation}
\delta_0((x_i)_{1\leq i\leq m})=m+\limsup_{\omega\rightarrow
0}\frac{\chi ((x_i+\omega s_i)_{1\leq i\leq m}:(s_i)_{1\leq i\leq
m})}{|\log\omega |},
\end{equation}
while its free entropy dimension is
\begin{equation}
\delta ((x_i)_{1\leq i\leq m})=m+\limsup_{\omega\rightarrow
0}\frac{\chi ((x_i+\omega s_i)_{1\leq i\leq m})}{|\log\omega |},
\end{equation}
where $(x_i)_{1\leq i\leq m}$ and the semicircular system
$(s_i)_{1\leq i\leq m}$ are free. Both free entropy dimensions can
be determined with the following formulae if the family
$(x_i)_{1\leq i\leq m}$ is free:
\begin{equation}
\delta_0((x_i)_{1\leq i\leq m})=\delta ((x_i)_{1\leq i\leq
m})=\sum_{i=1}^m\delta (x_i),
\end{equation}
\begin{equation}
\delta (x)=1-\sum_{s\in\mb{R}}(\mu (\{s\}))^2.
\end{equation}
We mention in this context the important Semicontinuity Problem
(\cite{13}): if $x_i$ is the $\mbox{SOT}$-limit of
$x_i^{(p)}\in\mc{P}$ as $p\rightarrow\infty$ (for all $1\leq i\leq
m$), does it follow then that
$\liminf_{p\rightarrow\infty}\delta_0((x_i^{(p)})_{1\leq i\leq
m})\geq\delta_0((x_i)_{1\leq i\leq m})$? An affirmative answer to
this question would imply (\cite{13}) the nonisomorphism of
$\mc{L}(\mb{F}_n)$ and $\mc{L}(\mb{F}_m)$ for $n\not= m$.

\section{Estimate of free entropy}\label{s2}

Lemma \ref{l1} gives an estimate for the free entropy of an
arbitrary system of generators of $\mc{M}$ which can be
$\omega$-approximated in the $2$-norm by certain noncommutative
polynomials. Typically, this situation is encountered under the
hypothesis of Theorem \ref{t1}, where the $\omega$-approximations
hold for all $\omega>0$.
\begin{lema}\label{l1}
Let $x_1,\ldots , x_m$ be self-adjoint generators of a
$\mbox{\!I\!I}_1$-factor $(\mc{M},\tau )$. Assume that $\mc{N}$ is
a subfactor of $\mc{M}$ and
$\mc{P}=\mc{P}_1\vee\mc{P}_2\subset\mc{N}$ is a von Neumann
subalgebra, where $\mc{P}_1=\mc{R}_1'\cap\mc{P}$,
$\mc{P}_2=\mc{R}_2'\cap\mc{P}$ and $\mc{R}_1$, $\mc{R}_2$ are
mutually commuting hyperfinite subfactors of $\mc{P}$. Assume
moreover that there exist self-adjoint elements
$m_j^{(e)},z_k\in\mc{M}$ (for $1\leq j\leq r+1$, $1\leq e\leq 2$,
$1\leq k\leq 2v$), mutually orthogonal projections $p_q\in\mc{M}$
(for $1\leq q\leq u$), projections $(p^{(t)})_t\subset\mc{P}_1$,
$(q^{(s)})_s\subset \mc{P}_2$ of trace $\frac{1}{2}$, and
noncommutative polynomials
$\Phi_{ji}^{(e)}((p^{(t)})_t,(q^{(s)})_s,(z_k)_k,(p_q)_q)$ which
are linear combinations of monomials of the form
$p^{(t_1)}q^{(s_1)}\ldots p^{(t_a)}q^{(s_a)}z_kp_q$, such that for
some $\omega>0$ and for all $1\leq i\leq m$
\begin{eqnarray}&&
\bigg|\bigg|\,
x_i-\frac{1}{2}\sum_{e=1}^{2}\sum_{j=1}^{r+1}\left(m_j^{(e)}
\Phi_{ji}^{(e)}\left((p^{(t)})_t,(q^{(s)})_s,(z_k)_k,(p_q)_q\right)\right.\\\nonumber&&
\hspace{1 cm}
\left.+\Phi_{ji}^{(e)}\left((p^{(t)})_t,(q^{(s)})_s,(z_k)_k,(p_q)_q\right)^*
m_j^{(e)}\right)\bigg|\bigg|_2<\omega\,.
\end{eqnarray}
Then
\begin{eqnarray}\label{1}&&
\chi\left((x_i)_{1\leq i\leq m}\right)=\chi\bigg((x_i)_{1\leq
i\leq m}:\left(m_j^{(e)} \right)_{j,e},(p^{(t)})_t,
(q^{(s)})_s,\\\nonumber&& \hspace{.5 cm}(z_k)_k,(p_q)_q\bigg)\leq
C(m,r,v,K)+(m-2r-2v-4)\log\omega\,.
\end{eqnarray}
where $C(m,r,v,K)$ is a constant depending only on $m$, $r$, $v$,
and $K=1+\max_{i,j,e}\left\lbrace\left|\left|
\Phi_{ji}^{(e)}\left((p^{(t)})_t,(q^{(s)})_s,(z_k)_k,(p_q)_q\right)\right|\right|_2,
||x_i||, \left|\left|m_j^{(e)}\right|\right|\right\rbrace$.
\end{lema}
\begin{proof} Let $c_0\geq 1$ be a fixed integer. Let
$\mc{M}_1\subset\mc{R}_1$ and $\mc{M}_2\subset\mc{R}_2$ be two
subalgebras isomorphic to $\mc{M}_{c_0}({\mb{C}})$ and let
$(e_{gh})_{g,h}$, $(f_{gh})_{g,h}$ be matrix units for $\mc{M}_1$
and $\mc{M}_2$ respectively. Consider a matricial microstate
$$\left((A_i)_i,(M_j^{(e)})_{j,e},(P^{(t)})_t,(Q^{(s)})_s,(Z_k)_k,(P_q)_q,
(E_{gh})_{g,h},(F_{gh})_{g,h}\right)$$ from the set of matricial
microstates
$$\Gamma_R\left((x_i)_i,(m_j^{(e)})_{j,e},(p^{(t)})_t,(q^{(s)})_s,
(z_k)_k,(p_q)_q,(e_{gh})_{g,h},(f_{gh})_{g,h};a,c,\epsilon\right).$$
We can assume (\cite{13}) that $||A_i||,||M_j^{(e)}||\leq K$. If
$a$ is large and $\epsilon >0$ is small enough, then $\forall
1\leq i\leq m$
\begin{eqnarray}&&
\bigg|\bigg|\,A_i-\frac{1}{2}\sum_{e=1}^{2}\sum_{j=1}^{r+1}\left(M_j^{(e)}
\Phi_{ji}^{(e)}\left((P^{(t)})_t,(Q^{(s)})_s,(Z_k)_k,(P_q)_q\right)\right.\\\nonumber
&&\hspace{1 cm}\left.
+\Phi_{ji}^{(e)}\left((P^{(t)})_t,(Q^{(s)})_s,(Z_k)_k,(P_q)_q\right)^*
M_j^{(e)}\right)\bigg|\bigg|_2<\omega
\end{eqnarray}
and
$\left|\left|\Phi_{ji}^{(e)}((P^{(t)})_t,(Q^{(s)})_s,(Z_k)_k,(P_q)_q)\right|\right|_2\leq
K$ for all $i,j,e$. For any $\delta >0$ there exists an injective
$*$-homomorphism
$\alpha_c:\mc{M}_1\vee\mc{M}_2\rightarrow\mc{M}_c({\mb{C}})$ such
that $$ \left|\left|\alpha_c(e_{gh})-E_{gh}\right|\right|_2<\delta
,\,
\left|\left|\alpha_c(e_{gh})-E_{gh}\right|\right|_2<\delta\,\forall
g,h$$ for large $a\in\mb{N}$ and small $\epsilon>0$, but
independently of $c$. The conditional expectation from $\mc{M}$
onto $\mc{M}_1'\cap\mc{M}$ is given by
$$E_{\mc{M}_1'\cap
\mc{M}}(x)=\frac{1}{c_0}\sum_{g,h=1}^{c_0}e_{gh}xe_{hg}=\Xi
(x,(e_{gh})_{g,h}).$$ Denote $P_0^{(t)}=\Xi
(P^{(t)},(\alpha_c(e_{gh}))_{g,h})\in
\alpha_c(\mc{M}_1)'\cap\mc{M}_c({\mb{ C}})$. For any $\delta_1>0$
and any $a_1\in\mb{N}$, since
$p^{(t)}=E_{\mc{M}_1'\cap\mc{M}}(p^{(t)})=\Xi
(p^{(t)},(e_{gh})_{g,h})$, it follows that
$$\left|\tau_c\left((P_0^{(t)})^l\right)-\tau \left((p^{(t)})^l\right)\right|<\delta_1
\forall 1\leq l\leq a_1$$ if $\epsilon ,\delta>0$ are small and
$a\in\mb{N}$ is large enough. Given $\delta_2>0$, if $\delta_1$ is
sufficiently small and $a_1$ is sufficiently large, there exists
(\cite{12}) a projection $P_1^{(t)}\in
\alpha_c(\mc{M}_1)'\cap\mc{M}_c({\mb{C}})$ of rank $\frac{c}{2}$,
such that
$\left|\left|P_1^{(t)}-P_0^{(t)}\right|\right|_2<\delta_2$. Note
that
$$\left|\left|P_0^{(t)}-P^{(t)}\right|\right|_2^2=\tau_c\left(\left(P^{(t)}-\Xi
(P^{(t)},(\alpha_c(e_{gh}))_{g,h})\right)^2\right)$$ and also
$$\tau \left(\left(p^{(t)}-\Xi
(p^{(t)},(e_{gh})_{g,h})\right)^2\right)=0,$$ hence
$\left|\left|P^{(t)}-P_0^{(t)}\right|\right|_2<\delta_2$ and thus
$\left|\left|P^{(t)}-P_1^{(t)}\right|\right|_2<2\delta_2$ if
$\epsilon ,\delta$ are small enough and $a$ is sufficiently large.
In this way we can find projections $(P_1^{(t)})_t\subset
\alpha_c(\mc{M}_1)'\cap \mc{M}_c({\mb{C}})$, $(Q_1^{(s)})_s\subset
\alpha_c(\mc{M}_2)'\cap\mc{M}_c({\mb{ C}})$, of rank
$\frac{c}{2}$, such that $ \left|\left|P^{(t)}-\right.\right.$
$\left.\left.P_1^{(t)} \right|\right|_2<2\delta_2$ and
$\left|\left|Q^{(s)}-Q_1^{(s)}\right|\right|_2<2\delta_2$ for all
$t,s$. If $\delta_2$ is small enough then we have moreover
$\forall 1\leq i\leq m$
\begin{eqnarray}&&
\bigg|\bigg|\,A_i-\frac{1}{2}\sum_{e=1}^{2}\sum_{j=1}^{r+1}\left(M_j^{(e)}
\Phi_{ji}^{(e)}\left((P_1^{(t)})_t,(Q_1^{(s)})_s,(Z_k)_k,(P_q)_q\right)
\right.\\\nonumber&&\hspace{1 cm}
\left.+\Phi_{ji}^{(e)}\left((P_1^{(t)})_t,(Q_1^{(s)})_s,(Z_k)_k,(P_q)_q\right)^*
M_j^{(e)}\right)\bigg|\bigg|_2<\omega .
\end{eqnarray}
Fix two copies $\mc{G}_1(c)\subset\alpha_c(\mc{M}_1)'\cap
\mc{M}_c({\mb{C}})$ and
$\mc{G}_2(c)\subset\alpha_c(\mc{M}_2)'\cap\mc{M}_c({\mb{ C}})$ of
the Grassmann manifold
$\mc{G}\left(\frac{c}{c_0},\frac{c}{2c_0}\right)$ and note that
there exists a unitary $U\in\mc{U}(c)$ such that
$UP_1^{(t)}U^*\in\mc{G}_1(c)$ and $UQ_1^{(s)}U^*\in\mc{G}_2(c)$
for all $t,s$. Lemma 4.3 in \cite{12} implies that given
$\delta_3>0$, there exist $a',c'\in\mb{N}$, $\epsilon_1>0$ such
that if $c\geq c'$ and if $(P_1,\ldots
,P_u)\in\Gamma_R((p_q)_q;a',c,\epsilon_1)$, then there exist
mutually orthogonal projections $P_1',\ldots ,P_u'\subset
\mc{M}_c^{sa}(\mb{C})$ such that $\mbox{rank}(P_q')=[\tau (p_q)c]$
and $\left|\left|P_q-P_q'\right|\right|_2<\delta_3$ $\forall 1\leq
q\leq u$. Let $(S_q)_q$ be fixed mutually orthogonal projections
with $\mbox{rank}(S_q)=[\tau (p_q)c]$ $\forall 1\leq q\leq u$ and
let $W\in\mc{U}(c)$ be a unitary such that $P_q'=W^*S_qW$ $\forall
1\leq q\leq u$. If $\delta_3>0$ is sufficiently small, then one
has
\begin{eqnarray}&&\bigg|\bigg|\,
UA_iW^*-\frac{1}{2}\sum\limits_{e=1}^{2}
\sum\limits_{j=1}^{r+1}\big(UM_j^{(e)}U^*
\Phi_{ji}^{(e)}\big((UP_1^{(t)}U^*)_t,(UQ_1^{(s)}U^*)_s,
\\\nonumber & &\hspace{.5 cm}(UZ_kW^*)_k,
(S_q)_q\big)+UW^*\Phi_{ji}^{(e)}\big((UP_1^{(t)}U^*)_t,
(UQ_1^{(s)}U^*)_s,\\\nonumber & &\hspace{.5 cm}
(UZ_kW^*)_k,(S_q)_q\big)^*
UM_j^{(e)}U^*(UW^*)\big)\bigg|\bigg|_2<\omega\,\forall 1\leq i\leq
m.
\end{eqnarray}
Consider a minimal $\theta$-net $(V_b)_{b\in B(c,K)}$ in $\{B\in
\mc{M}_c^{sa}({\mb{C}}): ||B||\leq K\}$ and a minimal
$\frac{\omega}{2K}$-net $(U_t)_{t\in T(c)}$ in $\mc{U}(c)$ with
respect to the uniform norm. Let also $(G_a^{(1)})_{a\in A(c)}$
and $(G_a^{(2)})_{a\in A(c)}$ be two minimal $\eta$-nets (relative
to the euclidian norm induced from $\mc{M}_c({\mb{C}})$) in
$\mc{G}_1(c)$ and respectively, $\mc{G}_2(c)$. From \cite{8} we
have $|T(c)|\leq (\frac{2CK}{\omega})^{c^2}$, $|B(c,K)|\leq
(\frac{CK}{\theta})^{c^2+c}$, $|A(c)|\leq
(\frac{C\sqrt{c}}{\eta})^{\frac{c^2}{2c_0^2}}$, where $C$ is a
universal constant. There exist indices $t,s\in T(c)$, $b(j,e)\in
B(c,K)$, $a(t),a(s)\in A(c)$ such that $\forall 1\leq i\leq m$
\begin{eqnarray}&&
\bigg|\bigg|\, U_tA_iW_s^*
-\frac{1}{2}\sum_{e=1}^{2}\sum_{j=1}^{r+1}\bigg(V_{b(j,e)}
\Phi_{ji}^{(e)}\bigg(\left(G_{a(t)}^{(1)}\right)_t,\left(G_{a(s)}^{(2)}\right)_s,
(T_k)_k,\\\nonumber&&\hspace{1 cm}
(S_q)_q\bigg)+U_tW_s^*\Phi_{ji}^{(e)}\bigg(\left(G_{a(t)}^{(1)}\right)_t,
\left(G_{a(s)}^{(2)}\right)_s, (T_k)_k,(S_q)_q\bigg)^*\\\nonumber
&&\hspace{1 cm}\cdot
V_{b(j,e)}U_tW_s^*\bigg)\bigg|\bigg|_e<\omega\sqrt{c}+\omega\sqrt{c}
+\frac{1}{2}\cdot 2(r+1)\cdot\left[\theta
K\sqrt{c}\right.\\\nonumber &&\hspace{1 cm}
\left.+D(\Phi)K\eta\sqrt{\alpha +\beta}+
2\cdot\frac{\omega}{2K}K^2\sqrt{c} +D(\Phi)K\eta\sqrt{\alpha
+\beta}\right.\\\nonumber&&\hspace{1 cm}
\left.+K\left(\theta\sqrt{c}+2K\frac{\omega}{2K}
\sqrt{c}\right)\right]
=2\left[K(r+1)+1\right]\omega\sqrt{c}\\\nonumber&&\hspace{1 cm}
+2\theta(r+1)K\sqrt{c} +2 D(\Phi)K\eta(r+1)\sqrt{\alpha +\beta}\,,
\end{eqnarray}
where $D(\Phi)$ is a Lipschitz constant depending on the $\Phi$'s,
and $\alpha ,\beta$ are the number of $P_1^{(t)}$'s and
$Q_1^{(s)}$'s. Choose $\theta =\frac{\omega}{2K(r+1)}$ and $\eta
=\frac{\omega\sqrt{c}}{2D(\Phi)K(r+1)\sqrt{\alpha +\beta}}$, so
that $|B(c,K)|\leq \left(\frac{2CK^2(r+1)}{\omega}\right)^{c^2+c}$
and $|A(c)|\leq$ $\left(\frac{2CD(\Phi)K(r+1)\sqrt{\alpha
+\beta}}{\omega}\right)^{\frac{c^2}{2c_0^2}}$. The volume of the
set of matricial microstates can be estimated as follows:
\begin{eqnarray}&&
\mbox{vol}_{mc^2}\big(\Gamma_R\big((x_i)_i:\left(m_j^{(e)}\right)_{j,e},(p^{(t)})_t,
(q^{(s)})_s,(z_k)_k,(p_q)_q,(e_{gh})_{g,h},\\\nonumber
&&\hspace{.5 cm} (f_{gh})_{g,h};a,c,\epsilon\big)\big)\leq
\left(\frac{2CK}{\omega}\right)^{2c^2}\cdot
\left(\frac{2CK^2(r+1)}{\omega} \right)^{2(r+1)(c^2+c)}\\\nonumber
&&\hspace{.5 cm}\cdot \left(\frac{2CD(\Phi)K(r+1)\sqrt{\alpha
+\beta}}{\omega}\right)^{\frac{c^2(\alpha +\beta )}{2c_0^2}}
\cdot\vol_{d_c} \left(0,(K+\mu)\sqrt{mc}\right)\\\nonumber
&&\hspace{.5 cm}\cdot\vol_{mc^2-d_c} \left(0, \mu\sqrt{mc}\right),
\end{eqnarray}
where $\mu=2\omega\left[K(r+1)+2\right]$ and $d_c$ denotes the
dimension of the range of the linear map that sends $(T_k)_k$ to
\begin{eqnarray}&& \bigg(\frac{1}{2}\sum\limits_{e=1}^{2}
\sum\limits_{j=1}^{r+1}\left(U_t^*V_{b(j,e)}
\Phi_{ji}^{(e)}\left(\left(G_{a(t)}^{(1)}\right)_t,\left(G_{a(s)}^{(2)}\right)_s,
(T_k)_k,(S_q)_q\right)W_s\right.\\\nonumber &&\hspace{.5 cm}
\left.+W_s^*\Phi_{ji}^{(e)}\left(\left(G_{a(t)}^{(1)}\right)_t,
\left(G_{a(s)}^{(2)}\right)_s, (T_k)_k,(S_q)_q\right)^*
V_{b(j,e)}U_t\right)\bigg)_{1\leq i\leq m}.
\end{eqnarray}
Since $(x_i)_{1\leq i\leq m}$ generates $\mc{M}$, the last
inequality implies the free entropy estimate
\begin{eqnarray}
\chi\left((x_i)_{1\leq i\leq m}\right)&=&\chi\big((x_i)_{1\leq
i\leq m}:\left(m_j^{(e)} \right)_{j,e},(p^{(t)})_t,
(q^{(s)})_s,\\\nonumber && (z_k)_k,(p_q)_q\big)\\\nonumber
&=&\chi\big((x_i)_{1\leq i\leq
m}:\left(m_j^{(e)}\right)_{j,e},(p^{(t)})_t,
(q^{(s)})_s,\\\nonumber && (z_k)_k,(p_q)_q,
(e_{gh})_{g,h},(f_{gh})_{g,h}\big)\\\nonumber &\leq&
C(m,r,v,K)+(m-2r-2v-4)\log\omega\,.
\end{eqnarray}
\end{proof}

\section{Applications}\label{s3}

The main application of the free entropy estimate from the
previous section is Theorem \ref{t1}: the (modified) free entropy
dimension of any set of generators of $\mc{M}$ is $\leq 2r+2v+4$
if the correspondence $_\mc{P}L^2(\mc{N},\tau )_\mc{A}$ is
$\mc{M}$-weakly contained in  a finitely generated
subcorrespondence $_\mc{P}H_\mc{A}$ of $_\mc{P}L^2(\mc{M},\tau
)_\mc{A}$, where $\mc{P}\subset\mc{N}$ is generated by the
(diffuse) relative commutants of two commuting copies of the
hyperfinite $\mbox{\!I\!I}_1$-factor $\mc{R}$,
$\mc{A}\subset\mc{N}$ is an abelian subalgebra, $v$ is the number
of vectors which span $_\mc{P}H_\mc{A}$, and $r$ is the integer
part of the Jones index $[\mc{M}:\mc{N}]$. The results concerning
the free group subfactors (such as the absence of regular
non-prime subfactors) are listed in Theorem \ref{t2}. We proceed
with two short Lemmas, \ref{l3} and \ref{l4}, which will be used
further in the proofs of Theorems \ref{t1} and \ref{t2},
respectively.
\begin{lema}\label{l3}
Let $_\mc{P}H_\mc{A}$ and $_\mc{P}K_\mc{A}$ be subcorrespondences
of $_\mc{P}L^2(\mc{M},\tau )_\mc{A}$ such that $_\mc{P}H_\mc{A}$
is generated by $v$ vectors and $_\mc{P} K_\mc{A}$ is
$\mc{M}$-weakly contained in $_\mc{P}H_\mc{A}$. Then $\forall
\epsilon>0$ $\forall \lambda_1,\ldots,\lambda_m\in K$ $\exists$
unitary $u\in\mc{M}$ $\exists \kappa_1,\ldots ,\kappa_v\in K$,
$\exists$ finite $\left\{b_{j,l}^{(i)}
\right\}_{i,j,l}\subset\mc{P}$, $\exists$ finite
$\left\{a_{j,l}^{(i)}\right\}_{i,j,l} \subset\mc{A}$ such that
$K=uKu^*$ and
$$\bigg|\bigg|u^*\lambda_iu-\sum_{j,l}b_{j,l}^{(i)}\kappa_la_{j,l}^{(i)}\bigg|\bigg|
<\epsilon\,\forall 1\leq i\leq m.$$
\end{lema}
\begin{proof} Note first that $_\mc{P}K_\mc{A}$ $\mc{M}$-weakly
contained in $_\mc{P}H_\mc{A}$ implies $\widetilde{_\mc{P}
K_\mc{A}}\in
V_{\mc{M}}\left(\widetilde{_\mc{P}H_\mc{A}},\epsilon_0,F,E,S\right)$
for all $\epsilon_0>0$ and all finite subsets $F\subset\mc{P}$,
$E\subset\mc{A}$, $S\subset H$. This shows in particular that
there exists a unitary $w\in\mc{M}$ such that $K=wHw^*$, hence
$\lambda_1=w\eta_1w^*,\ldots,\lambda_v=w\eta_vw^*$ for some
$\eta_1,\ldots,\eta_v\in H$. Since $_\mc{P}H_\mc{A}$ is generated
by $v$ vectors, there exist $\xi_1,\ldots, \xi_v\in H$ such that
$$_\mc{P}H_\mc{A}=_\mc{P}\overline{\mbox{sp}}^{||\cdot ||}
\left(\mc{P}\xi_1\mc{A}+\ldots +\mc{P}\xi_v\mc{A}\right)_\mc{A}.$$
Therefore, given $\epsilon_1>0$, there exist finite subsets
$F=\left\{b_{j,l}^{(i)}\right\}_{i,j,l} \subset\mc{P}$,
$E=\left\{a_{j,l}^{(i)}\right\}_{i,j,l}\subset\mc{A}$ such that
$$\bigg|
\bigg|\eta_i-\sum_{j,l}b_{j,l}^{(i)}\xi_la_{j,l}^{(i)}\bigg|\bigg|<\epsilon_1\,
\forall 1\leq i\leq m.$$ Given $\epsilon_2>0$, since
$\widetilde{_\mc{P} K_\mc{A}}\in
V_{\mc{M}}\left(\widetilde{_\mc{P}H_\mc{A}},\epsilon_2,F
,E,\{\xi_l\}_l\right)$, there exists a unitary $v\in\mc{M}$ such
that $K=vHv^*$ and
$$\left|\left|bv\xi_lv^*a-vb\xi_lav^*\right|\right|<\epsilon_2\,\forall b\in F\, \forall
a\in E\,\forall 1\leq l\leq v.$$ Let $\kappa_l=v\xi_lv^*$ $\forall
1\leq l\leq v$ and note that one has ($\forall 1\leq i\leq m$) the
following estimate:
\begin{eqnarray}\label{eq10}&&
\bigg|\bigg|vw^*\lambda_iwv^*-\sum_{j,l}b_{j,l}^{(i)}\kappa_la_{j,l}^{(i)}\bigg|\bigg|=
\bigg|\bigg|\eta_i-\sum_{j,l}v^*b_{j,l}^{(i)}v\xi_lv^*a_{j,l}^{(i)}v
\bigg|\bigg|\\\nonumber & &\hspace{1 cm}
\leq\bigg|\bigg|\eta_i-\sum_{j,l}b_{j,l}^{(i)}\xi_la_{j,l}^{(i)}\bigg|\bigg|+
\sum_{j,l}\bigg|\bigg|b_{j,l}^{(i)}\xi_la_{j,l}^{(i)}-v^*b_{j,l}^{(i)}v\xi_lv^*
a_{j,l}^{(i)}v\bigg|\bigg|\\\nonumber & &\hspace{1 cm}
<\epsilon_1+\sum_{j,l}\bigg|\bigg|v
b_{j,l}^{(i)}\xi_la_{j,l}^{(i)}v^*-b_{j,l}^{(i)}v\xi_lv^*a_{j,l}^{(i)}
\bigg|\bigg|.
\end{eqnarray}
The last term in (\ref{eq10}) is smaller than $\epsilon$ if
$\epsilon_1$ and $\epsilon_2$ are sufficiently small.
\end{proof}
\begin{lema}\label{l4}
Let $\mc{P}$ be a von Neumann algebra with a matrix unit
$(e_{ij})_{1\leq i,j\leq k}\subset\mc{P}$ and let also $\mc{A}$ be
an abelian algebra with a projection $q\in\mc{A}$. Assume that the
correspondence $_\mc{P}H_\mc{A}$ is finitely generated:
$H=\overline{\mbox{sp}}^{||\cdot||} \left(\mc{P}\xi_1\mc{A}+\ldots
+\mc{P}\xi_v\mc{A}\right)$ for some $\xi_1,\ldots,\xi_v\in H$.
Then the correspondence $_{\mc{P}_p}(pHq)_{\mc{A}_q}$ is also
finitely generated: $$pHq=\overline{\mbox{sp}}^{||\cdot||}
\left(\sum_{i=1}^k\sum_{l=1}^v\mc{P}_p\xi_{il}\mc{A}_q\right),$$
where $p=e_{11}$ and $\xi_{il}=e_{1i}\xi_lq$ $\forall 1\leq i\leq
k$ $\forall 1\leq l\leq v$.
\end{lema}
\begin{proof}
\begin{eqnarray}&&
pHq =\overline{\mbox{sp}}^{||\cdot||}
\left(\sum_{l=1}^vp\mc{P}\xi_l\mc{A}q\right)=
\overline{\mbox{sp}}^{||\cdot||}
\left(\sum_{i=1}^k\sum_{l=1}^vp\mc{P}e_{ii}\xi_l\mc{A}q\right)\\
\nonumber &&\hspace{1 cm}=\overline{\mbox{sp}}^{||\cdot||}
\left(\sum_{i=1}^k\sum_{l=1}^vp\mc{P}e_{i1}e_{1i}\xi_lq\mc{A}\right)
=\overline{\mbox{sp}}^{||\cdot||}
\left(\sum_{i=1}^k\sum_{l=1}^v\mc{P}_p\xi_{il}\mc{A}_q\right).
\end{eqnarray}
\end{proof}
\begin{teo}\label{t1}
Let $(\mc{M},\tau )$ be a $\mbox{\!I\!I}_1$-factor generated by
the self-adjoint elements $x_1,\ldots ,x_m$. Assume that
$\mc{N}\subset\mc{M}$ is a subfactor with the integer part of the
Jones index $[\mc{M}:\mc{N}]$ equal to $r$, $\mc{A}\subset\mc{N}$
is an abelian subalgebra, and $\mc{P}\subset\mc{N}$ is a
subalgebra such that $\mc{P}=\mc{P}_1\vee\mc{P}_2$, where
$\mc{P}_1=\mc{R}_1'\cap\mc{P}$, $\mc{P}_2=\mc{R}_2'\cap\mc{P}$ and
$\mc{R}_1,\mc{R}_2\subset\mc{P}$ are mutually commuting
hyperfinite subfactors. Assume moreover that $\mc{P}_1$,
$\mc{P}_2$ are diffuse von Neumann subalgebras and that the
correspondence $_\mc{P}L^2(\mc{N},\tau )_\mc{A}$ is
$\mc{M}$-weakly contained in  a subcorrespondence
$_\mc{P}H_\mc{A}$ of $_\mc{P}L^2(\mc{M},\tau )_\mc{A}$, generated
by $v$ vectors. Then
\begin{equation}
\delta_0(x_1,\ldots ,x_m)\leq 2r+2v+4\,.
\end{equation}
\end{teo}
\begin{proof} Note first that one has $\delta_0(x_1,\ldots ,x_m)\leq m$
(\cite{13}) and thus one can assume $m>2r+2v+4$. There exist
(\cite{6}) $m_1,\ldots ,m_{r+1} \in\mc{M}$ such that
\begin{equation}
x=\sum_{j=1}^{r+1}m_jE_\mc{N}(m_j^*x)\,\forall x\in\mc{M},
\end{equation}
where $E_\mc{N}:\mc{M}\rightarrow\mc{N}$ is the conditional
expectation onto $\mc{N}$. Use Lemma \ref{l3} to conclude that for
every $\epsilon>0$ there exist a unitary $u\in\mc{M}$,
self-adjoint vectors $\eta_1,\ldots,\eta_{2v}\in L^2(\mc{N},\tau
)^{sa}$ and finite subsets
$\left\{b_{p,k}^{(i,j)}\right\}_{i,j,p,k}\subset\mc{P}$,
$\left\{a_{p,k}^{(i,j)}\right\}_{i,j,p,k}\subset\mc{A}$ such that
$$\bigg|\bigg|u^*E_\mc{N}(m_j^*x_i)u-\sum_{k=1}^{2v}\sum_{p=1}^lb_{p,k}^{(i,j)}
\eta_ka_{p,k}^{(i,j)}\bigg|\bigg|_2<\epsilon\,\forall 1\leq i\leq
m\,\forall 1\leq j\leq r+1.$$ Since $u\mc{A}u^*$ is abelian, there
exist projections $p_1,\ldots ,p_u\in u\mc{A}u^*$ of sum $1$ such
that every $ua_{p,k}^{(i,j)}u^*$ is approximated sufficiently well
in the $||\cdot ||$-norm by linear combinations of these
projections. Being diffuse, both $u\mc{P}_1u^*$ and $u\mc{P}_2u^*$
are generated by their projections of trace $\frac{1}{2}$, hence
each $ub_{p,k}^{(i,j)}u^*$ is the SOT-limit of a sequence of
noncommutative polynomials
$\Psi_{p,k}^{(i,j)}\left((p^{(t)})_t,(q^{(s)})_s\right)$ in
projections of trace $\frac{1}{2}$, $(p^{(t)})_t\subset
u\mc{P}_1u^*$, $(q^{(s)})_s\subset u\mc{P}_2u^*$. Moreover,
$u\mc{N}^{sa}u^*$ is dense in $L^2(u\mc{N}u^*,\tau )^{sa}$ hence
there exist $z_1,\ldots ,z_{2v}$ self-adjoint elements of
$u\mc{N}u^*$ such that
\begin{equation}
\Lambda_{ji}\left((p^{(t)})_t,(q^{(s)})_s,(z_k)_k,(p_q)_q\right)
=\sum_{p=1}^l\sum_{k=1}^{2v}\sum_{q=1}^u\Psi_{p,k}^{(i,j)}\left((p^{(t)})_t,(q^{(s)})_s
\right)z_kp_q
\end{equation}
is sufficiently close to $E_\mc{N}(m_j^*x_i)$ in the $||\cdot
||_2$-norm. Therefore each $x_i$ can be approximated arbitrarily
well in the $||\cdot ||_2$-norm by elements of the form
\begin{equation}
\sum_{j=1}^{r+1}m_j\Lambda_{ji}\left((p^{(t)})_t,(q^{(s)})_s,(z_k)_k,(p_q)_q\right).
\end{equation}
Denote $m_j^{(1)}=\frac{m_j+m_j^*}{2}$,
$m_j^{(2)}=\frac{m_j-m_j^*}{2\sqrt{-1}}$, and
\begin{eqnarray}&&
\Phi_{ji}^{(1)}\left((p^{(t)})_t,(q^{(s)})_s,(z_k)_k,(p_q)_q\right)
=\Lambda_{ji}\left((p^{(t)})_t,(q^{(s)})_s,(z_k)_k,(p_q)_q\right)\\\nonumber
&&\hspace{1 cm}=-\sqrt{-1}
\Phi_{ji}^{(2)}\left((p^{(t)})_t,(q^{(s)})_s,(z_k)_k,(p_q)_q\right)\,.
\end{eqnarray}
Hence for every $\omega >0$ there exist polynomials
$\Phi_{ji}^{(e)} \big((p^{(t)})_t,(q^{(s)})_s,(z_k)_k,$ $
(p_q)_q\big)$ that are linear combinations of monomials of the
form $p^{(t_1)}q^{(s_1)}$ $\ldots p^{(t_a)}q^{(s_a)}z_kp_q$ such
that
\begin{eqnarray}\label{2}&&
\bigg|\bigg|\,x_i-\frac{1}{2}\sum_{e=1}^2\sum_{j=1}^{r+1}\left(m_j^{(e)}\Phi_{ji}^{(e)}
\left((p^{(t)})_t,(q^{(s)})_s,(z_k)_k,(p_q)_q\right)\right.\\
\nonumber &&\hspace{.5 cm}+\left.\Phi_{ji}^{(e)}
\left((p^{(t)})_t,(q^{(s)})_s,(z_k)_k,(p_q)_q\right)^*m_j^{(e)}\right)\bigg|\bigg|_2
<\omega\,\forall 1\leq i\leq m\,.
\end{eqnarray}
Also, one can assume that $\left(\left|\left|\Phi_{ji}^{(e)}
\left((p^{(t)})_t,(q^{(s)})_s,(z_k)_k,
(p_q)_q\right)\right|\right|_2\right)_{i,j,e}$ are uniformly
bounded by a constant $D$ depending only on
$(||m_j^*x_i||)_{i,j}$. Consider a semicircular system
$(s_i)_{1\leq i \leq m}$, free from $(x_i)_{1\leq i\leq m}$. Note
that since $\left(m_j^{(e)}\right)_{j,e}$, $(p^{(t)})_t$,
$(q^{(s)})_s$, $(z_k)_k$, $(p_q)_q$ are all contained in
$\{x_i+\omega s_i, s_i:1\leq i\leq m\}''$, one has (\cite{13})
\begin{eqnarray}\label{3}&&
\chi\left((x_i+\omega s_i)_{1\leq i\leq m}:(s_i)_{1\leq i\leq
m}\right) =\chi\bigg((x_i+\omega s_i)_{1\leq i\leq m}:(s_i)_{1\leq
i\leq m},\\\nonumber &&\hspace{.5 cm}\left( m_j^{(e)}
\right)_{j,e},(p^{(t)})_t, (q^{(s)})_s, (z_k)_k, (p_q)_q\bigg)\leq
\chi\bigg((x_i+\omega s_i)_{1\leq i\leq m}:\\\nonumber
&&\hspace{.5 cm}\left(
m_j^{(e)}\right)_{j,e},(p^{(t)})_t,(q^{(s)})_s, (z_k)_k,
(p_q)_q\bigg)\,\forall 1\leq i\leq m.
\end{eqnarray}
The inequalities (\ref{2}) imply
\begin{eqnarray}&&
\bigg|\bigg|x_i+\omega
s_i-\frac{1}{2}\sum_{e=1}^2\sum_{j=1}^{r+1}\left(m_j^{(e)}
\Phi_{ji}^{(e)}
\left((p^{(t)})_t,(q^{(s)})_s,(z_k)_k,(p_q)_q\right)\right.\\\nonumber
&&\hspace{.5 cm} +\left.\Phi_{ji}^{(e)}
\left((p^{(t)})_t,(q^{(s)})_s,(z_k)_k,(p_q)_q\right)^*m_j^{(e)}\right)\bigg|
\bigg|_2<2\omega\,\forall 1\leq i\leq m
\end{eqnarray}
hence, by (\ref{3}) and the free entropy estimate (\ref{1}) from
Lemma \ref{l1},
$$
\chi\left((x_i+\omega s_i)_{1\leq i\leq m}:(s_i)_{1\leq i\leq
m}\right) \leq C(m,r,v,K)+(m-2r-2v-4)\log 2\omega\,.
$$
The estimate for the (modified) free entropy dimension follows now
immediately:
\begin{eqnarray}&&
\delta_0(x_1,\ldots,x_m)=m+\limsup_{\omega\rightarrow
0}\frac{\chi\left((x_i+ \omega s_i)_{1\leq i\leq m}:(s_i)_{1\leq
i\leq m}\right)}{|\log\omega |}\\\nonumber &&\hspace{1 cm}\leq m
+\limsup_{\omega\rightarrow 0}\frac{C(m,r,v,K)+(m-2r-2v-4)\log
2\omega}{| \log\omega |}\\\nonumber &&\hspace{1 cm}=2r+2v+4\,.
\end{eqnarray}
\end{proof}
\begin{coro}\label{c2}
Let $(\mc{M},\tau )$ be a $\mbox{\!I\!I}_1$-factor generated by
the self-adjoint elements $x_1,\ldots ,x_m$. Assume that
$\mc{N}\subset\mc{M}$ is a subfactor with the integer part of the
Jones index $[\mc{M}:\mc{N}]$ equal to $r$, $\mc{A}\subset\mc{N}$
is an abelian subalgebra, and $\mc{P}\subset\mc{N}$ is a
subalgebra such that $\mc{P}=\mc{P}_1\vee\mc{P}_2$, where
$\mc{P}_1=\mc{R}_1'\cap\mc{P}$, $\mc{P}_2=\mc{R}_2'\cap\mc{P}$ and
$\mc{R}_1,\mc{R}_2\subset\mc{P}$ are mutually commuting
hyperfinite subfactors. Assume moreover that $\mc{P}_1$,
$\mc{P}_2$ are diffuse von Neumann subalgebras and that
$L^2(\mc{N},\tau
)=\overline{\mbox{sp}}^{||\cdot||_2}(\mc{P}\xi_1\mc{A}+ \ldots
+\mc{P}\xi_v\mc{A})$ for some vectors $\xi_1,\ldots ,\xi_v\in
L^2(\mc{N},\tau )$. Then
\begin{equation}
\delta_0(x_1,\ldots ,x_m)\leq 2r+2v+4\,.
\end{equation}
\end{coro}

L. Ge and S. Popa proved (\cite{4}) that if $\mc{Q}$ is a finite
von Neumann algebra with no atoms and with a faithful normal trace
$\tau:\mc{Q}\rightarrow\mb{C}$ and if, moreover,
$\alpha:\Gamma\rightarrow\mbox{Aut}(\mc{Q})$ is a trace-preserving
properly outer action of a countable discrete group $\Gamma$ on
$\mc{Q}$, then there exist an abelian subalgebra
$\mc{A}\subset\mc{Q}$ and $\xi\in
L^2(\mc{Q}\times_\alpha\Gamma,\tau)$ such that
$L^2(\mc{Q}\times_\alpha\Gamma,
\tau)=\overline{\mbox{sp}}^{||\cdot||_2}\mc{Q}\xi\mc{A}$. In the
same vein, one has the following Lemma:
\begin{lema}\label{l2}
Let $\mc{Q}=\mc{Q}_1\vee\mc{Q}_2\simeq\mc{Q}_1\otimes\mc{Q}_2$ be
a non-prime subfactor of a $\mbox{\!I\!I}_1$-factor $\mc{N}$. If
$\mc{Q}$ is regular in $\mc{N}$, then there exist diffuse abelian
subalgebras $\mc{A}_1\subset\mc{Q}_1$, $\mc{A}_2\subset\mc{Q}_2$
and an abelian subalgebra $\mc{A}_3\subset\mc{Q}'\cap\mc{N}$ such
that the correspondence $_{\mc{P}}L^2(\mc{N},\tau)_\mc{A}$ is
cyclic i.e., $_{\mc{P}}L^2(\mc{N},
\tau)_\mc{A}=\overline{\mbox{sp}}^{||\cdot||_2}\mc{P}\xi\mc{A}$
for some $\xi\in L^2(\mc{N},\tau)$, where
$\mc{P}=\mc{Q}\vee(\mc{Q}'\cap\mc{N})$ and
$\mc{A}=\mc{A}_1\vee\mc{A}_2\vee\mc{A}_3$.
\end{lema}
\begin{proof} Let $\Gamma$ be a countable group of unitaries in
$N_\mc{N}(\mc{Q})$ such that
$\overline{\mbox{sp}}^{||\cdot||_2}\Gamma =L^2(\mc{N},\tau)$. Note
that
$\mc{P}\simeq\mc{Q}_1\otimes\mc{Q}_2\otimes(\mc{Q}'\cap\mc{N})$,
$N_\mc{N}(\mc{Q})\subset N_\mc{N}(\mc{P})$ and
$\mc{P}'\cap\mc{N}\subset\mc{P}$. Use $\S2$ in \cite{14} to
conclude that there exist maximal abelian subalgebras
$\mc{A}_1\subset\mc{Q}_1$, $\mc{A}_2\subset\mc{Q}_2$,
$\mc{A}_3\subset\mc{Q}'\cap\mc{N}$ with the property that for
every finite subset $W\subset\mbox{sp}\Gamma$ ($=$linear span of
$\Gamma$) and every $\epsilon>0$ there exists a finite partition
of the identity with projections $(p_i)_{i\in I}
\subset\mc{A}=\mc{A}_1\vee\mc{A}_2\vee\mc{A}_3$ such that
$||\sum_{i\in I}p_iwp_i-E_\mc{A}(w)||_2<\epsilon$ $\forall w\in
W$, where $E_\mc{A}:\mc{N}\rightarrow\mc{A}$ is the conditional
expectation onto $\mc{A}$. Pick a vector $\xi\in L^2(\mc{N},\tau)$
such that $(\xi,u) \not= 0$ $\forall u\in\Gamma$. In fact, one can
assume that $\xi =v\in\mc{U} (\mc{N})$ since the set of unitaries
$v\in\mc{U}(\mc{N})$ such that $(v,u)=\tau(vu^*) \not=0$ $\forall
u\in\Gamma$ is a $G_\delta$-dense subset in $\mc{U}(\mc{N})$
(\cite{4}). Let $(w_n)_{n\geq 1}\subset\mbox{sp}\Gamma$ be an
orthonormal basis of $L^2(\mc{N},\tau)$ and write $\xi
=\sum_{n\geq 1}\alpha_nw_n$, $\xi_m =\sum_{n\geq 1}^m\alpha_nw_n$,
where $\alpha_n=(\xi,w_n)\in\mb{C}$ $\forall n\geq 1$. Given
$\delta
>0$, there exists $m\geq 1$ such that $||\xi-\xi_m||_2<\delta$.
For $u\in\Gamma$ and $\epsilon >0$ let $(p_i)_{i\in I}$ be a
finite partition of the identity with projections from $\mc{A}$
such that $||\sum_{i\in
I}up_iu^*w_np_i-uE_\mc{A}(u^*w_n)||_2<\epsilon$ $\forall 1\leq
n\leq m$. Note that
\begin{eqnarray}&&
\bigg|\bigg|\sum_{i\in I}up_iu^*\xi
p_i-uE_\mc{A}(u^*\xi)\bigg|\bigg|_2\leq \bigg|\bigg|\sum_{i\in
I}up_iu^*(\xi -\xi_m) p_i\bigg|\bigg|_2\\\nonumber &&\hspace{1
cm}+ \bigg|\bigg|\sum_{i\in I}up_iu^*\xi_m p_i-uE_\mc{A}(u^*\xi_m)
\bigg|\bigg|_2+\bigg|\bigg|uE_\mc{A}(u^*(\xi
-\xi_m))\bigg|\bigg|_2\\\nonumber &&\hspace{1 cm}
<\delta+\sum_{n=1}^m|\alpha_n|\epsilon+\delta ,
\end{eqnarray}
hence $u\in\overline{\mbox{sp}}^{||\cdot||_2}\mc{P}\xi\mc{A}$
since $\epsilon, \delta>0$ can be chosen arbitrarily small,
$up_iu^*\in\mc{P}$, $p_i\in\mc{A}$ $\forall i\in I$, $\mc{A}$ has
no atoms, and $\tau(\xi u^*)\not=0$.
\end{proof}
\begin{teo}\label{t2}
Let $\mc{N}$ be a subfactor of finite index in the interpolated
free group factor $\mc{M}=\mc{L}(\mb{F}_t)$ ($1<t\leq\infty$) and
let also $r$ denote the integer part of the index. The following
statements are true: \newline i) $\mc{N}$ does not have regular
non-prime subfactors; \newline ii) the correspondence
$_\mc{P}L^2(\mc{N},\tau)_\mc{A}$ is not finitely generated if
$\mc{A}$ is an abelian subalgebra of $\mc{N}$ and
$\mc{P}=\mc{P}_1\vee\mc{P}_2$ is a subalgebra of $\mc{N}$ such
that $\mc{P}_1=\mc{R}_1'\cap\mc{P}$ and
$\mc{P}_2=\mc{R}_2'\cap\mc{P}$ are both diffuse, $\mc{R}_1$, $
\mc{R}_2$ are mutually commuting hyperfinite subfactors of
$\mc{P}$, and $\mc{R}_1\cap\mc{A}$, $\mc{R}_2\cap\mc{A}$ have
projections of arbitrarily small trace;
\newline
iii) the correspondence $_\mc{P}L^2(\mc{N},\tau)_\mc{A}$ is not
$\mc{M}$-weakly contained in any finitely generated correspondence
$_\mc{P}H_\mc{A}$ if $2r+2v+4<t\leq\infty$, $\mc{A}$ is an abelian
subalgebra of $\mc{N}$ and $\mc{P}=\mc{P}_1\vee\mc{P}_2$ is a
subalgebra of $\mc{N}$ such that $\mc{P}_1=\mc{R}_1'\cap\mc{P}$
and $\mc{P}_2=\mc{R}_2'\cap\mc{P}$ are both diffuse and
$\mc{R}_1$, $ \mc{R}_2$ are mutually commuting hyperfinite
subfactors of $\mc{P}$;
\newline iv) $\mc{N}$ does not have regular diffuse hyperfinite
$*$-subalgebras (DHSA).
\end{teo}
\begin{proof} i) Assume that $\mc{N}$ has a regular nonprime subfactor
$\mc{Q}=\mc{Q}_1\vee\mc{Q}_2\simeq\mc{Q}_1\otimes\mc{Q}_2$ and
denote $\mc{P}=\mc{Q}\vee (\mc{Q}'\cap\mc{N})$. By Lemma \ref{l2},
there exist diffuse abelian subalgebras $\mc{A}_1\subset\mc{Q}_1$,
$\mc{A}_2\subset\mc{Q}_2$, an abelian subalgebra
$\mc{A}_3\subset\mc{Q}'\cap\mc{N}$, and a vector $\xi\in
L^2(\mc{N},\tau)$ such that
$L^2(\mc{N},\tau)=\overline{\mbox{sp}}^{||\cdot
||_2}\mc{P}\xi\mc{A}$, where
$\mc{A}=\mc{A}_1\vee\mc{A}_2\vee\mc{A}_3$. We consider first the
case $\mc{M}=\mc{L}(\mb{F}_t)$ with $1<t<\infty$. Since $\mc{A}_1$
and $\mc{A}_2$ are diffuse, for any $k\geq 1$, there exist
projections $p_1\in\mc{A}_1$, $p_2\in\mc{A}_2$ such that $\tau
(p_1)=\tau (p_2)=\frac{1}{k}$. Let $(e_{ij})_{1\leq i,j\leq
k}\subset\mc{Q}_1$, $(f_{ls})_{1\leq l,s\leq k}\subset\mc{Q}_2$ be
two matrix units such that $e_{11}=p_1$, $f_{11}=p_2$ and denote
$p=p_1p_2$. Use Lemma \ref{l4} to conclude
$$
L^2(\mc{N}_p,\tau_p)=pL^2(\mc{N},\tau)p=\overline{\mbox{sp}}^{||\cdot
||_2}\sum_{1\leq i,l\leq k}\mc{P}_p\xi_{il}\mc{A}_p\,,
$$
where $\xi_{il}=e_{1i}f_{1l}\xi p$ $\forall 1\leq i,l\leq k$. Let
$\mc{R}_1\subset (\mc{Q}_2)_p$, $\mc{R}_2\subset (\mc{Q}_1)_p$ be
hyperfinite subfactors and denote $\mc{P}_1=\mc{R}_1'\cap\mc{P}_p$
and $\mc{P}_2=\mc{R}_2'\cap\mc{P}_p$. If the integer part of
$[\mc{M}:\mc{N}]$ is equal to $r$, then the integer part of
$[\mc{M}_p:\mc{N}_p]$ is also equal to $r$ and the estimate of
free entropy dimension from Corollary \ref{c2} implies
\begin{equation}\label{5}
\delta_0(x_1,\ldots,x_m)\leq 2r+2k^2+4
\end{equation}
for any system $(x_1,\ldots ,x_m)$ of self-adjoint generators of
$\mc{M}_p$. On the other hand, by the compression formula
(\cite{1}, \cite{7}),
$$\mc{M}_p\simeq\mc{L}\left(\mb{F}_{1+(t-1)\tau(p)^{-2}}\right)=
\mc{L}\left(\mb{F}_{1+(t-1)k^4}\right).$$ Moreover (\cite{12},
\cite{13}), $\mc{L}\left(\mb{F}_{1+(t-1)k^4}\right)$ has a system
of generators $(x_1,$ $\ldots ,x_m)$ with $\delta_0(x_1,\ldots
,x_m)=1+(t-1)k^4$, hence the inequality (\ref{5}) implies
$1+(t-1)k^4\leq 2r+2k^2+4$ which is of course impossible if $k$ is
sufficiently large.

Let us consider now the case $\mc{M}=\mc{L}(\mb{F}_\infty)$, when
(\cite{11}) $\mc{M}$ is generated by an infinite semicircular
system $(x_i)_{i\geq 1}$. With the estimate of free entropy
(\ref{1}) we conclude that there exist elements $\left(m_j^{(e)}
\right)_{j,e}$, $(p^{(t)})_t$, $(q^{(s)})_s$, $(z_k)_k,(p_q)_q$
(as stated in the proof of Theorem \ref{t1}) such that
\begin{eqnarray}\nonumber
\chi\left((x_i)_{1\leq i\leq m}:\left(m_j^{(e)}
\right)_{j,e},(p^{(t)})_t,
(q^{(s)})_s,(z_k)_k,(p_q)_q\right)<\chi\left((x_i)_{1\leq i\leq
m}\right).
\end{eqnarray}
Let $\mc{M}_n=\{x_1,\ldots ,x_n\}''$ and
$E_n:\mc{M}\rightarrow\mc{M}_n$ be the conditional expectation
onto $\mc{M}_n$. Since
\begin{eqnarray}&&
\bigg((x_i)_{1\leq i\leq m},\left(E_n\left(m_j^{(e)}\right)
\right)_{j,e},\left(E_n\left(p^{(t)}\right)\right)_t,
\left(E_n\left(q^{(s)}\right)\right)_s,\\\nonumber &&\hspace{1 cm}
\left(E_n\left(z_k\right) \right)_k,
\left(E_n\left(p_q\right)\right)_q\bigg)_{n\geq 1}
\end{eqnarray}
converges in distribution as $n\rightarrow\infty$ to
$$\left((x_i)_{1\leq i\leq m},\left(m_j^{(e)}
\right)_{j,e},(p^{(t)})_t, (q^{(s)})_s,(z_k)_k,(p_q)_q\right)\,$$
there exists an integer $n>m$ such that
\begin{eqnarray}&&
\chi\bigg((x_i)_{1\leq i\leq m}:\left(E_n\left(m_j^{(e)}\right)
\right)_{j,e},\left(E_n\left(p^{(t)}\right)\right)_t,
\left(E_n\left(q^{(s)}\right)\right)_s,\\\nonumber &&\hspace{1 cm}
\left(E_n\left(z_k\right)
\right)_k,\left(E_n\left(p_q\right)\right)_q\bigg)
<\chi\left((x_i)_{1\leq i\leq m}\right),
\end{eqnarray}
hence
\begin{eqnarray}&&
\chi\left((x_i)_{1\leq i\leq n}\right)=\chi\bigg((x_i)_{1\leq
i\leq n}:\left(E_n\left(m_j^{(e)}\right)
\right)_{j,e},\left(E_n\left(p^{(t)}\right)\right)_t,\\\nonumber&&\hspace{1
cm}
\left(E_n\left(q^{(s)}\right)\right)_s,\left(E_n\left(z_k\right)
\right)_k,\left(E_n\left(p_q\right)\right)_q\bigg)\\\nonumber
&&\hspace{1 cm}\leq\chi\bigg((x_i)_{1\leq i\leq
m}:\left(E_n\left(m_j^{(e)}\right)
\right)_{j,e},\left(E_n\left(p^{(t)}\right)\right)_t,\\\nonumber&&\hspace{1
cm}
\left(E_n\left(q^{(s)}\right)\right)_s,\left(E_n\left(z_k\right)
\right)_k,\left(E_n\left(p_q\right)\right)_q\bigg)\\\nonumber
&&\hspace{1 cm}+ \chi\left(x_{m+1},\ldots,x_n\right)
<\chi\left(x_1,\ldots,x_m\right)+\chi\left(x_{m+1},\ldots,x_n\right),
\end{eqnarray}
contradiction.
\newline
ii) The statement is a direct consequence of the free entropy
dimension estimate from Corollary \ref{c2} if $2r+2v+4<t<\infty$.
If $1<t\leq 2r+2v+4$, first cut down by a projection $p=p_1p_2$
with $p_1\in\mc{R}_1\cap\mc{A}$, $p_2\in\mc{R}_1\cap\mc{A}$ of
sufficiently small trace. Note that this increases $\delta_0$ as
in the proof of i) and use then Lemma \ref{l4} and the free
entropy dimension estimate from Corollary \ref{c2}. The case
$t=\infty$ can be treated as in the proof of i).
\newline
iii) The case $2r+2v+4<t<\infty$ is consequence of the estimate of
free entropy dimension from Theorem \ref{t1}. The case $t=\infty$
can be also treated as in the proof of i).
\newline
iv) Let $\mc{Q}\subset\mc{N}$ be a regular DHSA of $\mc{N}$. As in
the proof of Lemma \ref{l2}, conclude that there exist a diffuse
abelian subalgebra $\mc{A}_1\subset\mc{Q}$ and an abelian
subalgebra $\mc{A}_3\subset \mc{Q}'\cap\mc{N}$ such that
$_\mc{P}L^2(\mc{N},\tau)_\mc{A}= \overline{\mbox{sp}}^{||\cdot
||_2}\mc{P}\xi\mc{A}$ for some $\xi\in L^2(\mc{N},\tau)$, where
$\mc{P}=\mc{Q}\vee(\mc{Q}'\cap\mc{N})$ and
$\mc{A}=\mc{A}_1\vee\mc{A}_3$. Since $\mc{Q}$ is a DHSA of
$\mc{N}$, this implies (with the notations from Theorem \ref{t1})
that for any $\epsilon>0$ there exist mutually commuting
hyperfinite subfactors $\mc{R}_1,\mc{R}_2\subset\mc{N}$ (depending
on $\epsilon$) such that
$$\mbox{dist}_{||\cdot ||_2}\left(E_\mc{N}(m_j^*x_i),\mbox{sp}\mc{R}_1\mc{R}_2\xi\mc{A}\right)
<\epsilon\,\forall 1\leq i\leq m\,\forall 1\leq j\leq r+1.$$ As in
the proof of Theorem \ref{t1}, one obtains the free entropy
dimension estimate $\delta_0(x_1,\ldots,x_m)\leq 2r+6$ and then
iv) follows from this estimate in a fashion similar to the proof
of i).
\end{proof}
\begin{coro}\label{c1}
The subfactors $\mc{N}$ of finite index in the interpolated free
group factors $\mc{L}(\mb{F}_t)$ ($1<t\leq\infty$) are not crossed
products of nonprime subfactors or diffuse hyperfinite subalgebras
by properly outer actions of countable discrete groups.
\end{coro}
\begin{proof} Let $\mc{Q}\subset\mc{N}$ be either a nonprime subfactor or
a diffuse hyperfinite subalgebra. Recall that if $\Gamma$ is a
countable discrete group and if
$\alpha:\Gamma\rightarrow\mbox{Aut}(\mc{Q})$ is a properly outer
action of $\Gamma$ on $\mc{Q}$ such that
$\mc{N}\simeq\mc{Q}\times_\alpha\Gamma$  then $\mc{Q}$ is regular
in $\mc{N}$. Use then i) and iv) from Theorem \ref{t2}.
\end{proof}

\end{document}